\documentclass[a4paper,reqno,11pt,final]{amsart}
\usepackage{latexsym}
\usepackage{listings}
\usepackage{amsthm,amsmath,amssymb}
\usepackage[dvips]{graphicx, color}
\usepackage[title]{appendix}
\usepackage{float}
\usepackage[linesnumbered,ruled,vlined]{algorithm2e}
\usepackage{color}
\usepackage{bbm}

\usepackage{tikz}
\usetikzlibrary{quantikz}
\usepackage{mathtools}

\newtheorem{thm}{Theorem}

\newtheorem{mdef}{Definition}

\newtheorem{prop}[thm]{Proposition}

\newtheorem{cor}[thm]{Corollary}
\newtheorem{conj}[thm]{Conjecture}
\newtheorem{lem}[thm]{Lemma}

\newtheorem{rem}[]{Remark}

\newcommand{\cE}{{\mathcal E}}
\newcommand{\cF}{{\mathcal F}}

\newcommand{\cL}{{\mathcal L}}
\newcommand{\cM}{{\mathcal M}}

\newcommand{\cO}{{\mathcal O}}

\newcommand{\cV}{{\mathcal V}}

\newcommand{\C}{{\Bbb C}}

\newcommand{\F}{{\Bbb F}}

\newcommand{\Q}{{\Bbb Q}\hspace{.06em}}
\newcommand{\R}{{\Bbb R}}

\newcommand{\Z}{{\Bbb Z}}

\newcommand{\frm}{{\frak m}}

\def\={\:=\:}  \def\+{\,+\,}

\def\a{\alpha} \def\b{{\beta}}  \def\ba{\overline\a}   \def\z{\zeta}

\def\be{\begin{equation}}   \def\ee{\end{equation}}
\def\bes{\begin{equation*}}   \def\ees{\end{equation*}}
\def\ba{\begin{aligned}}   \def\ea{\end{aligned}}
\def\bc{\begin{cases}}   \def\ec{\end{cases}}
\def\bp{\begin{proof}}   \def\ep{\end{proof}}

\newcommand{\Res}{\mathrm{Res}}

\newcommand{\Aut}{\mathrm{Aut}}

\def\SL{\mathrm{SL}}

\def\qqan{\qquad\mathrm{and}\qquad}
\def\qan{\quad\mathrm{and}\quad}

\def\SL{\mathrm{SL}}

\def\Lra{\Longrightarrow}

\def\bbm1{\mathbbm 1}

\def\wh{\widehat}
\def\be{\begin{equation}}   \def\ee{\end{equation}}
\def\bes{\begin{equation*}}   \def\ees{\end{equation*}}
\def\bea{\begin{equation}\begin{aligned}}   
\def\eea{\end{aligned}\end{equation}}

\def\whv{\wh{v}}

\def\Pic{\mathrm{Pic}}

\def\om{\omega}

\def\bm{\begin{matrix}}
\def\em{\end{matrix}}
\def\bpm{\begin{pmatrix}}
\def\epm{\end{pmatrix}}
\def\Hom{\mathrm{Hom}}

\def\rk{\mathrm{rank}}

\def\bl{\big(}
\def\br{\big)}

\begin{document}
\title[Riemann Hypothesis for Non-Abelian Zetas of Genus 2 Curves]{\bf Riemann Hypothesis for Non-Abelian Zeta Functions of Genus 2 Curves} 
\author{\bf Zhan SHI}  
\address{}
\date{}
\maketitle
\begin{abstract}
In this paper, we investigate Weng zeta functions associated with curves of genus $2$ over finite fields. Building upon Weng’s framework for non-abelian zeta functions, we establish that, as the rank $n$ tends to infinity, the Riemann Hypothesis holds for these zeta functions. Our proof relies on the geometric properties of the moduli space of semi-stable bundles, together with several established results for high rank zeta functions, complemented by detailed asymptotic analysis. This result provides new evidence supporting the general validity of the Riemann Hypothesis for Weng zeta functions and offers insight into the analytic structure of non-abelian zeta functions associated with higher-genus curves.
\end{abstract}



\section{Introduction}
The study of zeta functions lies at the heart of modern number theory and arithmetic geometry. Since Riemann’s seminal work on the classical zeta function for number fields, the distribution of its nontrivial zeros—the Riemann Hypothesis—has been regarded as one of the most profound unsolved problems in mathematics. Since the 20th century, the zeta function has been generalized to function fields, aiming to capture deeper arithmetic and geometric structures. 

Weng’s framework of high rank zeta functions provides a rich and unified theory that generalized classical Artin zeta functions for curves and introduced new families of zeta functions over function fields. These non-abelian functions satisfy analogues of the functional equation and demonstrate connections between geometric and arithmetic structures. An outstanding question in this theory is whether such Weng zeta functions satisfy Riemann Hypothesis.

In this paper, we investigate Weng zeta functions associated with curves of genus $2$ over finite fields. Building upon Weng’s framework for non-abelian zeta functions, we establish that, as the rank $n$ tends to infinity, the Riemann Hypothesis holds for these zeta functions. Our proof relies on the geometric properties of the moduli space of semi-stable bundles, together with several established results for high rank zeta functions, complemented by detailed asymptotic analysis. This result provides new evidence supporting the general validity of the Riemann Hypothesis for Weng zeta functions and offers insight into the asymptotically analytic structure of non-abelian zeta functions associated with higher-genus curves.

The structure of this paper is as follows.
\begin{itemize}
\item Section 2 recalls the definition of Weng zeta functions and reviews some basic and advanced properties necessary for our study. It also reviews the current statements of Riemann Hypothesis for zeta function for curves over finite fields that has been proven.
\item Section 3 discusses the properties of $\a$- and $\b$- invariants, including the counting miracle and the semi-stable mass, emphasizing their roles in our study.
\item Section 4 focuses on the case of curves of genus 2, fixing $g=2$ to simplify the general framework and make the analysis more explicit.
\item Section 5 establishes the main theorem, showing that the Riemann Hypothesis holds for genus 2 Weng zeta functions in the large-rank limit.
\end{itemize}

\section{Non-Abelian Zeta Function for Function Fields}
\subsection{Artin Zeta Function}
Let $X$ be a smooth projective curve of genus $g$ over a finite field $\F_q$ with $q$ elements. If $P$ is a closed point of $X/\F_q$ and $\cO_P$(resp. $\frm_P$) is its valuation ring(resp. maximal ideal), then we define its degree by
$$
d(P)=[k_P:\F_q],
$$
where $k_P$ denotes the residue class field $\cO_P/\frm_P$. As usual, the norm of $P$ is defined by
$$
N(P)=q^{d(P)}.
$$
For a divisor $D=\sum_Pn_PP$, we put
$$
d(D)=\sum_Pn_Pd(P) \qqan N(D)=q^{d(D)}.
$$
\begin{mdef}[Artin Zeta Function]
Let $X$ be a smooth projective curve of genus $g$ over $\F_q$. Then its Artin zeta function is defined as the following
\be\label{AZ}
\zeta_{X/\F_q}(s):=\sum_{D\ge 0}\frac{1}{N(D)^s},\quad \Re (s) >1,
\ee
where $D$ runs over all effective divisors of $X/\F_q$. Here an effective(or equivalently, positive) divisor means all $n_P\in\Z_{\ge 0}$.
\end{mdef}
As usual, making a variable transformation $t=q^{-s}$, we have 
$$
Z_{X/\F_q}(t):=\zeta_{X/\F_q}(s)=\sum_{D\ge 0}t^{d(D)}.
$$

A well-known result of Euler, the Euler product formula, has an analogue in the setting of function fields, which yields the following equation:
\be\label{EPF}
\sum_{D\ge 0}N(D)^{-s}=\prod_P(1-N(P)^{-s})^{-1}.
\ee
This is analogous to the definition of the classical Riemann zeta function $\z(s)$, since every effective divisor $D\ge0$ can be uniquely written as a finite sum of closed points $P$ on $X/\F_q$ with nonnegative integer coefficients, reflecting the unique factorization of divisors into closed points, or equivalently, prime divisors.

The following theorem is a famous result of Weil:
\begin{thm}[Zeta Properties]\label{AZP}
Let $X$ be a smooth projective curve of genus $g$ over a finite filed $\F_q$, and let $Z_{X/\F_q}(t)$ be its Artin zeta function. 
\begin{enumerate}
\item (Rationality)
$Z_{X/\F_q}(t)$ is a rational function of the form
$$
Z_{X/\F_q}(t)=\frac{P_{X/\F_q}(t)}{(1-t)(1-qt)},
$$
where $P_{X/\F_q}(t)$ is a polynomial of degree $2g$ with integer coefficients of the form
$$
P_{X/\F_q}(t)=1+a_1t+\dots+a_gt^g+a_{g-1}qt^{g+1}+\dots+q^gt^{2g}.
$$
\item (Functional Equation) 
The functional equation holds for $\zeta_{X/\F_q}(s)$
$$
\zeta_{X/\F_q}(1-s)=q^{(g-1)(2s-1)}\cdot \zeta_{X/\F_q}(s).
$$
\item (Riemann Hypothesis) 
The polynomial $P_{X/\F_q}(t)$ can be factored as 
$$
P_{X/\F_q}(t)=\prod_{i=1}^g(1-\om_{i,+}t)(1-\om_{i,-}t),
$$
while for all $i=1,2,\dots,g$, $|\om_{i,+}|=|\om_{i,+}|=\sqrt q$.
\end{enumerate}
\end{thm}

\subsection{Weng Zeta Function}
Non-abelian zeta function, introduced by Weng, is a generalization of Artin zeta function. To understand the definition, we first do the following transformation:
$$
\ba
\zeta_{X/\F_q}(t)&=\sum_{D\ge 0}t^{d(D)}
\\&=\sum_{[D]}\frac{\#\{D\in[D]:D\ \rm{positive}\}}{q-1}t^{d([D])}
\\&=\sum_{\cL\in\Pic(X/\F_q)}\frac{q^{h^0(X,\cL)}-1}{\#\Aut(\cL)}t^{d(\cL)},
\ea
$$
where $[D]$ denotes the equivalent class of divisor $D$, $\cL\in\Pic(X/\F_q)$ denotes the line bundle in the Picard group of $X/\F_q$, and $h^0(X,\cL)$ is the dimension of the 0-th cohomology group $H^0(X,\cL)$ of $\cL$ over $X$. The transformation works since we have the relation
$$
\Pic(X/\F_q)\cong Cl(X/\F_q).$$

This means that the classical Artin zeta function can be regarded as counting line bundles of $X/\F$. Thus, to generalize Artin zeta function, we may replace line bundles by rank $n$ vector bundles.

\begin{mdef}[Rank $n$ Zeta Function, see \cite{HRZ1} and \cite{HRZ2}]
For a fixed $n\in\Z_{\ge 1}$, the non-abelian zeta function of rank $n$ for $X$ over $\F_q$ is defined by 
\be\label{WZ}
\zeta_{X/\F_q,n}(s) := \sum_{\cV}\frac{q^{h^0(X,\cV)}-1}{\#Aut(\cV)}(q^{-s})^{d(\cV)},
\ee
where $\cV$ runs over all semi-stable vector bundles on $X/\F_q$ with rank $n$ and degree $d(\cV)\in n\Z_{\ge 0}$. 
\end{mdef}

We should notice that besides the formal summation, two additional conditions are included in the definition. The first condition, that is, $\cV$ runs over only semi-stable vector bundles of rank $n$ on $X/\F_q$, should be introduced, since the space of rank $n$ vector bundles with a fixed degree is unbounded and thus the summation doesn't work. The geometric stability condition was first introduced by Mumford, which is based on the (Mumford's) $\mu$-slope, defined by 
$$
\mu(\cV)=\frac{\deg{\cV}}{\rk{\cV}}.
$$
$\cV$ is called semi-stable if for all sub-bundle $\cV'$ of $\cV$, we have $\mu(\cV')\le\mu(\cV)$.

In fact, in Weng's first attempt, the second condition, that is, degree $d(\cV)$ is a multiple of $n$, hadn't been considered, even if most of zeta properties could be proved without this condition, but with the help of the standard zeta techniques and several algebraic geometric theorems, say Riemann--Roch theorem, duality theorem and vanishing theorem for semi-stable vector bundles. After considering some simple examples, Weng found that such defined function do not satisfy the Riemann hypothesis. This suggests that there may be something fundamentally wrong in that definition. Several years later, motivated by Drinfeld's paper on counting rank two irreducible representations of the fundamental groups of $X/\F_q$, the second condition was added by Weng. The current definition looks nice and works quite well.

\subsection{Properties of Rank $n$ Zetas}
To study the properties of rank $n$ zetas, we first introduce the following two fundamental non-abelian geo-arithmetic invariants of $X/\F_p$ associated to rank $n$ semi-stable vector bundles of degree $d$, which are defined by
$$
\ba
\alpha_{X/\F_q,n}(d)&:=\sum_{\cV\in\cM_{X/\F_q,n}(d)}\frac{q^{h^0(X,\cV)}-1}{\#\Aut(\cV)},
\\
\beta_{X/\F_q,n}(d)&:=\sum_{\cV\in\cM_{X/\F_q,n}(d)}\frac{1}{\#\Aut(V)},
\ea
$$
where $\cV$ runs over all semi-stable vector bundles of rank $n$ and degree $d$.
Then the zeta function $\zeta_{X/\F_q,n}(s)$ can be regarded as a generating function of $\alpha_{X/\F_q,n}(mn)$:
\be\label{GF}
\zeta_{X,n}(s)=\sum_{n|d}\alpha_{X/\F_q,n}(d)t^d=\sum_{m=0}^{\infty}\alpha_{X/\F_q,n}(mn)T^m,
\ee
where we set $Q=q^n$ and $T=t^n=Q^{-s}$.

The following theorem has been proved, using Riemann--Roch theorem, duality theorem and vanishing for semi-stable bundles:
\begin{thm}[High Rank Zeta Properties, see \cite{HRZ2}]\label{HRZP}
Let $X/\F_q$ be a projective regular curve of genus $g$ over finite field $\F_q$. Then the rank $n$ zeta function $\zeta_{X/\F_q,n}(s)$ satisfies the following zeta properties:
\begin{enumerate}
\item (Naturality) 
When $n=1$, $\zeta_{X/\F_q,1}(s)$ is actually the classical Artin zeta function $\zeta_{X/\F_q}(s)$ of $X/\F_q$.
\item (Rationality)
There exists a polynomial $P_{X/\F_q,n}(T)\in \Q[T]$ of degree $2g$, such that
$$
\zeta_{X/\F_q,n}(s)=Z_{X/\F_q,n}(T)=\frac{P_{X/\F_q,n}(T)}{(1-T)(1-QT)}.
$$
\item (Functional Equation)
The functional equation holds for $\wh\zeta_{X/\F_q,n}(s)$
$$
\wh\zeta_{X/\F_q,n}(1-s)=\wh\zeta_{X/\F_q,n}(s).
$$
where $\wh\zeta_{X/\F_q,n}(s):=T^{(1-g)}\cdot \zeta_{X/\F_q,n}(s)$.
\item (Residue in Geometry)
$Z_{X/\F_q,n}(T)$ admits only two simple poles at $T=1$ and $T=1/Q$, while their residues are given by
$$
\Res_{T=1}Z_{X/\F_q,n}(T)=-\Res_{T=1/Q}Z_{X/\F_q,n}(T)=\beta_{X/\F_q,n}(0).
$$
\end{enumerate}
\end{thm}

According to Riemann--Roch theorem and vanishing theorem for semi-stable vector bundles, we may obtain that when $d<0$, $\a_{X/\F_q,n}(d)=0$ and 
$$
\a_{X/\F_q,n}(d)=(q^{d-n(g-1)}-1)\b_{X/\F_q,n}(d),\ \ \forall d\geq 2n(g-1).
$$
In addition, directly from the definition, we have
$$
\b_{X/\F_q,n}(mn)=\b_{X/\F_q,n}(0),\ \ \forall m\in\Z.
$$

Therefore, form (\ref{GF}), we have the following
\begin{thm}[\cite{HRZ2}]\label{RoHRZ}
The rank $n$ zeta function $\wh\zeta_{X/\F_q,n}(s)$ of a genus $g$ projective regular integral curve $X$ over $\F_q$ is given by
$$
\ba
&\wh Z_{X/\F_q,n}(T):=\wh\zeta_{X/\F_q,n}(s)
\\=&\sum_{m=0}^{g-2}\a_{X/\F_q,n}(mn)\left(T^{m-(g-1)}+Q^{(g-1)-m}T^{(g-1)-m}\right) +\a_{X/\F_q,n}\bl n(g-1)\br
\\&+\frac{(Q-1)\b_{X/\F_q,n}(0)\cdot T}{(1-T)(1-QT)}.
\ea
$$
\end{thm} 

In particular, the degree $2g$ polynomial $P_{X/\F_q,n}(T)$ can be obtained from 
$$
\wh Z_{X/\F_q,n}(T)=\frac{P_{X/\F_q,n}(T)}{(1-T)(1-QT)\cdot T^{g-1}},
$$
that is,
\bea\label{PT}
&P_{X/\F_q,n}(T)=\wh Z_{X/\F_q,n}(T)\cdot(1-T)(1-QT)\cdot T^{g-1}
\\=&\Big(\sum_{m=0}^{g-2}\a_{X/\F_q,n}(nm)(T^m+Q^{g-1-m}T^{2g-2-m})+\a_{X/\F_q,n}(n(g-1))T^{g-1}\Big)
\\&\times(1-T)(1-QT)+(Q-1)T^g\b_{X/\F_q,n}(0).
\eea
One can easily check that $P_{X/\F_q,n}(T)$ is a degree $2g$ polynomial with rational coefficients whose leading term and constant term are $\a_{X/\F_q,n}(0)Q^gT^{2g}$ and $\a_{X/\F_q,n}(0)$ respectively.

\subsection{Riemann Hypothesis}
Riemann hypothesis is one of the most important zeta properties. The first result of Riemann hypothesis for zetas of curves over finite fields is the following Hasse's theorem on elliptic curves, first conjectured by Artin and proved by Hasse.

\begin{thm}[Hasse's Theorem on Elliptic Curves]
Let $E/\F_q$ be an elliptic curve over finite field $\F_q$, then we have the following estimates
$$
|\#E(\F_q)-(q+1)|\le2\sqrt{q},
$$
where $\#E(\F_q)$ denotes the number of points on the elliptic curve $E/\F_q$.
\end{thm}

This theorem can be regarded as the if and only if condition to the Riemann hypothesis for Artin zeta function $\zeta_{E/\F_q}(s)$, based on a discussion on the multiplicative structure of Artin zetas. This theorem had been generalized to higher genus algebraic curves $X/\F_q$ by Weil, which is often called Hasse-Weil bound. Accordingly, as introduced in Theorem \ref{AZP}, Riemann hypothesis for (rank 1) zeta function $\zeta_{X/\F_q}(s)$ is regarded as being proven.

\begin{thm}[Hasse-Weil Bound]\label{HWB}
Let $X/\F_q$ be a projective regular curve of genus $g$ over finite field $\F_q$, then we have the following estimates
$$
|\#X(\F_q)-(q+1)|\le2g\sqrt{q},
$$
where $\#X(\F_q)$ denotes the number of points on the curve $X/\F_q$.
\end{thm}

While for the high rank zeta function $\zeta_{X/\F_q,n}(s)$ for curves over finite fields, it is also conjectured that Riemann hypothesis holds.
\begin{conj}[Rank $n$ Riemann Hypothesis, see \cite{HRZ2} and \cite{RH}]
Let $X$ be an integral regular projective curve of genus $g$ over $\F_q$. Then all zeros of the rank n non-abelian zeta function $\zeta_{X/\F_q,n}(s)$ of $X/\F_q$ lie on the central line $\Re(s)=\frac{1}{2}$. That is, 
$$
\zeta_{X/\F_q,n}(s)=0\Lra \Re(s)=\frac{1}{2}.
$$
\end{conj}
Obviously, this is equivalent to the condition that the norm of reciprocal roots of $P_{X/\F_q,n}(T)$ are all $\sqrt{Q}$. In fact, according to definition and functional equation, we may write $P_{X/\F_q,n}(T)$ as 
\bea
P_{X/\F_q,n}(T)&=\a_{X/\F_q,n}(0)\prod_{i=1}^g(1-\omega_{X/\F_q,n,i}T)(1-\omega_{X/\F_q,n,i}^{\prime}T)
\\&=\a_{X/\F_q,n}(0)\prod_{i=1}^g(1-a_{X/\F_q,n,i}T+QT^2),
\eea
where
$$
\bc
\omega_{X/\F_q,n,i}\cdot\omega_{X/\F_q,n,i}^{\prime}=Q,
\\
\omega_{X/\F_q,n,i}+\omega_{X/\F_q,n,i}^{\prime}=a_{X/\F_q,n,i}.
\ec
$$

Thus, we have the following fundamental lemma.
\begin{lem}\label{iffRH}
Let $X$ be an integral  regular projective curve of genus $g$ over $\F_q$. Then the following conditions are equivalent:
\begin{enumerate}
\item [(1)] The rank $n$ Riemann Hypothesis holds for $\zeta_{X/\F_q,n}(s)$.
\item [(2)] For all $i=1,\dots g$, 
$$
|\omega_{X/\F_q,n,i}|=|\omega_{X/\F_q,n,i}^{\prime}|=\sqrt{Q}.
$$
\item[(3)] For all $i=1,\dots g$,
$$
\overline{\omega_{X/\F_q,n,i}}=\omega_{X/\F_q,n,i}^{\prime}.
$$
\item[(4)] For all $i=1,\dots g$,
$$
a_{X/\F_q,n,i}\in \R \qan a_{X/\F_q,n,i}\in(-2\sqrt{Q},2\sqrt{Q}).
$$
\item[(5)] For all $i=1,\dots g$,
$$
\omega_{X/\F_q,n,i}\in \C -\R \qan \omega_{X/\F_q,n,i}^{\prime}\in\C -\R.
$$
\end{enumerate}
\end{lem}

Although, as mentioned in Theorem\ \ref{HWB}, Riemann hypothesis for Artin zeta $\zeta_{X/\F_q}(s)$ had been proved, this still remains as a conjecture and wildly open for general cases of high rank zeta $\zeta_{X/\F_q,n}(s)$. The first breakthrough for high rank Riemann hypothesis is the following theorem.

\begin{thm}[Weng--Zagier's Theorem on Rank $n$ Riemann Hypothesis for Elliptic Curves, see \cite{EC}]\label{ECRH}
Let $E$ be an integral regular projective curve of genus $g=1$, that is, an elliptic curve, over finite field $\F_q$. Then rank $n$ Riemann hypothesis holds for $\zeta_{E/\F_q,n}(s)$.
\end{thm}
This theorem comes from a much stronger estimate of the linear coefficient in $P_{E/\F_q,n}(T)$, which indicates that $P_{E/\F_q,n}(T)$ admits no real zeros. The result relies on the fundamental properties of Atiyah bundles on elliptic curves, combined with a careful and detailed analysis of the underlying combinatorial structures. 

Furthermore, for $\zeta_{X/\F_q,n}(s)$ with $X$ of arbitrary genus $g$, the following cases for rank $2$ and $3$ have also been proved. Rank $2$ Riemann Hypothesis was proved by Yoshida, while Rank $3$ Riemann Hypothesis was proved by Weng. For more details, one may refer to \cite{R2RH} and \cite{RH}.

\begin{thm}[Rank $2\ \&$ Rank $3$ Riemann Hypothesis, see \cite{R2RH} $\&$ \cite{RH}]\label{R23RH}
Let $X$ be an integral regular projective curve of genus $g$ over finite field $\F_q$. Then $\zeta_{E/\F_q,n}(s)$ satisfies the rank $n$ Riemann hypothesis if $n\in\{2,\ 3\}$.
\end{thm}

\subsection{Special Uniformity of Zeta Functions}
The Special Uniformity theorem claims that for a global field $F$, the geometrically defined $\wh\zeta_{F,n}$ coincides with the Lie theoretically defined $\wh\zeta_F^{\SL_n}$. This has been verified for function fields in \cite{SU}, which enables us to decompose $\wh\zeta_{X/\F_q,n}$ into $\wh\zeta_{X/\F_q}$. 

We note that the proof of the Special Uniformity theorem is based on the equivalence between the formula of $\wh\zeta_{X/\F_q}^{\SL_n}(s)$ established in \cite{SU} and the formula of $\wh\zeta_{X/\F_q,n}(s)$, as given by of Mozgovoy and Reineke(Theorem 7.2 in \cite{MR}).
\begin{thm}[Special Uniformity, see Theorem 1, \cite{SU}]\label{SU}
For an integral  regular projective curve $X$ of genus $g$ over $\F_q$, we have
$$
\ba
\wh\zeta_{X/\F_q,n}(s)&=\wh\zeta_{X/\F_q}^{\SL_n}(s)
\\&=\sum_{a=1}^{n}\Biggl(\sum_{\substack{k_1,\ldots,k_p>0\\ k_1+\ldots+k_p=n-a}}\frac{\wh v_{k_1}\ldots\wh v_{k_p}}{\prod_{j=1}^{p-1}(1-q^{k_j+k_{j+1}})} \frac{1}{(1-q^{ns-n+a+k_{p}})}\Biggr)
\\&\times q^{\binom{n}{2}(g-1)}\wh\zeta_{X/\F_q}(ns-n+a)
\\&\times\Biggl(\sum_{\substack{l_1,\ldots,l_r>0\\ l_1+\ldots+l_r=a-1}}\frac{1}{(1-q^{-ns+n-a+1+l_{1}})}\frac{\wh v_{l_1}\ldots\wh v_{l_r}}{\prod_{j=1}^{r-1}(1-q^{l_j+l_{j+1}})}\Biggr),
\ea
$$
where we have defined
$$
\wh v_n:=\wh\zeta_{X/\F_q}^*(1)\wh\zeta_{X/\F_q}(2)\dots\wh\zeta_{X/\F_q}(n), \quad
\wh\zeta_{X/\F_q}^*(1):=\lim_{s\to1}(1-q^{1-s})\wh\zeta_{X/\F_q}(1).
$$
\end{thm}

To understand this formula, we first emphasize again that this formula enable us to decompose $\wh\zeta_{X/\F_q,n}$ into $\wh\zeta_{X/\F_q}$. In fact, if we define
$$
\ba
\wh\z_{X/\F_q,n}^{[a]}(s)
&:=\Biggl(\sum_{\substack{k_1,\ldots,k_p>0\\ k_1+\ldots+k_p=n-a}}\frac{\wh v_{k_1}\ldots\wh v_{k_p}}{\prod_{j=1}^{p-1}(1-q^{k_j+k_{j+1}})} \frac{1}{(1-q^{ns-n+a+k_{p}})}\Biggr)
\\&\times\wh\zeta_{X/\F_q}(ns-n+a)
\\&\times\Biggl(\sum_{\substack{l_1,\ldots,l_r>0\\ l_1+\ldots+l_r=a-1}}\frac{1}{(1-q^{-ns+n-a+1+l_{1}})}\frac{\wh v_{l_1}\ldots\wh v_{l_r}}{\prod_{j=1}^{r-1}(1-q^{l_j+l_{j+1}})}\Biggr),
\ea
$$ 
then $\wh\zeta_{X/\F_q,n}(s)$ can be simply written as
$$
\wh\zeta_{X/\F_q,n}(s)=q^{\binom{n}{2}(g-1)}\sum_{a=1}^{n}\wh\z_{X/\F_q,n}^{[a]}(s).
$$
Furthermore, for each $\wh\z_{X/\F_q,n}^{[a]}(s)$, it is a product of three functions. The middle one is the rank $1$ zeta function $\wh\zeta_{X/\F_q}(s)$ where $s=ns-n+a$. The rest parts are summations of rational functions which are related to ordered partition of $n-a$ and $a-1$, respectively. We should notice that ordered partition is a combination of integer partition and permutation. For example, for the first summation, the set of partition is $\{(k_1,\ldots,k_p)\}$ such that $k_i>0$ and $k_1+\ldots+k_p=n-a$. This means we have to fix the length $p$ of the partition at first, while certainly $p$ is from $1$ to $n-a$. Furthermore, since this is ordered partition, permutation of each partition $(k_1,\ldots,k_p)$ should be considered. Since $a$ is from $1$ to $n$, it may occur that the ordered partition is performed on 0. In this case, it will be regarded as $1$. Additionally, we mention that if the length $p=1$, then the term $\prod_{j=1}^{p-1}(1-q^{k_j+k_{j+1}})$ will be $1$ and the following $k_p$ will be $n-a$. 

\section{Properties of $\a$- and $\b$- Invariants}
As introduced in previous section, $\a$- and $\b$- invariants play a key roll in the theory of non-abelian zeta functions, since they can represent the coefficients of $\wh\z_{X/\F_q,n}$. These invariants have been extensively investigated, with a long-standing history and a well-developed theoretical framework(see e.g. \cite{EC}, \cite{SU}, \cite{HST}, \cite{SGHR} and \cite{MR}). In fact, for a general curve $X/\F_q$, $\beta_{X/\F_q,n}$ was first introduced by Harder-Narasimhan in \cite{HN}. Some basic properties, based on Riemann--Roch theorem and vanishing theorem, has been introduced before.  In this section, we will introduce some more properties of $\a$- and $\beta$- invariants in degree $0$, namely, $\a_{X/\F_q,n}(0)$ and $\b_{X/\F_q,n}(0)$, for our limited purpose. 

Recall that the definitions are
$$
\alpha_{X/\F_q,n}(d):=\sum_{\cV}\frac{q^{h^0(X,\cV)}-1}{\#\Aut(\cV)}
\qqan
\beta_{X/\F_q,n}(d)=\sum_{\cV}\frac{1}{\#\Aut(V)}.
$$

\subsection{Counting Miracle}
The first property we are going to introduce is the so-called {\it counting miracle}, which gives an intrinsic relation between $\a_{X/\F_q,n+1}(0)$ and $\b_{X/\F_q,n}(0)$. This was first conjectured by Weng, for the case $X$ is an elliptic curve, or equivalently, $g=1$, based on some careful computations for lower rank, with the help of the classification Atiyah bundles. Later, this was generalized for general (integral
regular) projective curves $X/\F_q$ of genus $g$ by Sugahara and independently
by Mozgovoy-Reineke. The key point is that the category of semi-stable bundles of degree 0 is an abelian category.
\begin{thm}[Counting Miracle, see \cite{SGHR}, and \cite{EC} for elliptic curve case]\label{CM}
Let $\cE_0$ be a stable vector bundle of rank $m(<n)$ and degree $0$. Then we have 
$$
\sum_{\cE\in\cM_{X/\F_{q,n}(0)}}\frac{q^{\#\Hom(\cE_0,\cE)}-1}{\#\Aut(\cE)}=q^{m(n-m)(g-1)}\sum_{\cF\in\cM_{X/\F_q,n-m}(0)}\frac{1}{\#\Aut(\cF)}.
$$
In particular, if we set $n=n+1$ and $m=1$, then we obtain
$$
\a_{X/\F_q,n+1}(0)=q^{n(g-1)}\beta_{X/\F_q,n}(0).
$$
\end{thm}

\subsection{Semi-Stable Mass}
In this subsection, in which $X/\F_q$ is again a curve of arbitrary genus $g\ge1$, we introduce a closed formula for $\beta_{X/\F_q,n}(0)$ for all $n\ge1$, which is called Semi-Stable Mass, combining results of \cite{HN}, \cite{PP}, and \cite{VF}. For more details, one may refer to \cite{EC}.

According to vanishing theorem, the invariant $\beta_{X/\F_q,n}(d)$ is periodic in $d$ of period $n$. We give a renormalized definition by setting
$$
\wh\beta_{X/\F_q,n}(d):=q^{-(g-1)n(n-1)/2}\beta_{X/\F_q,n}(d).
$$

\begin{thm}[Semi-Stable Mass, see \cite{EC}]\label{SSM}
For $n\ge1$ and any $d\in\Z$, we have
$$
\wh\beta_{X/\F_q,n}(d)=\sum_{k\ge1}(-1)^{k-1}\sum_{n_1+\dots+n_k=n}\prod_{j=1}^k\wh v_{n_j}\times\prod_{j=1}^{k-1}\frac{q^{(n_j+n_{j+1})\{d(n_1+\dots+n_j)/n\}}}{q^{n_j+n_{j+1}}-1},
$$
where we have defined
$$
\wh v_n:=\wh\zeta_{X/\F_q}^*(1)\wh\zeta_{X/\F_q}(2)\dots\wh\zeta_{X/\F_q}(n), \quad
\wh\zeta_{X/\F_q}^*(1):=\lim_{s\to1}(1-q^{1-s})\wh\zeta_{X/\F_q}(1),
$$
and $\{x\}$ denote the fractional part of $x$.

In particular, when $d=0$, we have 
$$
\wh\beta_{X/\F_q,n}(0)=\sum_{\substack{n_1,\ldots,n_k>0\\ n_1+\ldots+n_k=n}}\frac{\prod_{j=1}^k\wh v_{n_j}}{\prod_{j=1}^{k-1}(1-q^{n_j+n_{j+1}})}.
$$
\end{thm}

\begin{rem}
One may easily observe that the closed formula of $\wh\b_{X/\F_q,n}(0)$ also has the ordered partition structure, similar with what appears in Theorem\,\ref{SU}. In fact, they share the same algorithm.
\end{rem}

Combining Theorem\,\ref{CM} and Theorem\,\ref{SSM}, we could obtain closed formulae of both $\a$- and $\b$- invariants. However,  it depends on the integer partition and permutation, making it extremely complicated for a large $n$. Nevertheless, there is a beautiful recursion formula for $\beta_{E/\F_q,n}(0)$, proved in \cite{EC}, for the case $g=1$, or equivalently, the curve is an elliptic curve. The recursion formula for $\wh\beta_{X/\F_q,n}(0)$ of general curve hasn't been discovered.

Even if it's hard to obtain the exact values of $\a$- and $\b$- invariants, it should be more simple to consider the approximate values, or more general, the asymptotic behavior of these invariants, when $n$ tends to infinity. The asymptotic behavior should only depend on the curve $X/\F_q$ and the rank $n$. 

\subsection{Estimates of $\a$- and $\b$- Invariants}
In this subsection, we study the asymptotic behavior of $\a$- and $\b$- invariants in degree $0$ and give some basic estimates.

We start with the estimate of $\wh v_n$.
\begin{prop}[Estimate of $\wh v_n$]\label{Eovh}
Under the same notation introduced before, we have 
$$
\wh v_n=O(q^{(g-1)n(n+1)/2}).
$$
\end{prop}
\bp
Recall Theorem\,\ref{AZP}(1), we have
$$
\zeta_{X/\F_q}(s)=Z_{X/\F_q}(t)=\frac{P_{X/\F_q}(t)}{(1-t)(1-qt)},
$$
where $P_{X/\F_q}(t)$ is a polynomial of the form
$$
P_{X/\F_q}(t)=1+a_1t+\dots+a_gt^g+a_{g-1}qt^{g+1}+\dots+q^gt^{2g}.
$$
Then by definition, we obtain
$$
\whv_1=\wh\z_{X/\F_q}^*(1)=\frac{q^{g-1}P_{X/\F_q}(q^{-1})}{1-q^{-1}}.
$$
This means
$$
\whv_1=\wh\z_{X/\F_q}^*(1)=O(q^{g-1}).
$$
For $n\ge2$, we have
$$
\wh\z_{X/\F_q}(n)=\frac{q^{n(g-1)}P_{X/\F_q}(q^{-n})}{(1-q^{-n})(1-q^{1-n})}.
$$
Since $P_{X/\F_q}(q^{-n})$ is $O(1)$, we have
$$
\wh\z_{X/\F_q}(n)=O(q^{n(g-1)}).
$$
Finally, based on the definition, we could summarize that
$$
\wh v_n=O(q^{(g-1)+2(g-1)+...+n(g-1)})=O(q^{(g-1)n(n+1)/2})
$$
as wanted.
\ep

\begin{prop}[Estimate of $\a_{X/\F_q,n}(0)$, $\b_{X/\F_q,n}(0)$ and $\wh\b_{X/\F_q,n}(0)$]
For a genus $g$ curve $X/\F_q$ and a fixed rank $n$, we have: 
$$
\wh\b_{X/\F_q,n}(0)=O\Big(q^{(g-1)n(n+1)/2}\Big),
$$
$$
\b_{X/\F_q,n}(0)=q^{(g-1)\binom{n}{2}}\wh\b_{X/\F_q,n}(0)=O\Big(q^{(g-1)n^2}\Big).
$$
Then by counting miracle, we have
$$
\a_{X/\F_q,n}(0)=O\Big(q^{(g-1)(n^2-n)}\Big).
$$
\end{prop}
\bp
According to Theorem\,\ref{SSM}, we have
$$
\wh\beta_{X/\F_q,n}(0)=\sum_{\substack{n_1,\ldots,n_k>0\\ n_1+\ldots+n_k=n}}(-1)^{k-1}\frac{\prod_{j=1}^k\wh v_{n_j}}{\prod_{j=1}^{k-1}(q^{n_j+n_{j+1}}-1)}.
$$
Consider the numerator, we have
$$
\ba
\prod_{j=1}^k\wh v_{n_j}&=\prod_{j=1}^kO(q^{(g-1)n_j(n_j+1)/2})
\\&=O(q^{\sum_{j=1}^k(g-1)n_j(n_j+1)/2})
\\&=O(q^{\frac{(g-1)}{2}(n_1^2+n_2^2+...+n_k^2+n)})
\ea
$$
since we have $n_1+n_2+...+n_k=n$.

As for the denominator, we have
$$
\ba
\prod_{j=1}^{k-1}(q^{n_j+n_{j+1}}-1)&=\prod_{j=1}^{k-1}O(q^{n_j+n_{j+1}})
\\&=O(q^{\sum_{j=1}^{k-1}n_j+n_{j+1}})
\\&=O(q^{2n-n_1-n_k}).
\ea
$$

However, recall that if the length $k$ of the partition $(n_1,n_2,...,n_k)$ is $1$,  then the numerator is just $\whv_n=O(q^{(g-1)(n^2+n)/2})$, while the denominator is $1$.
In fact, since $n_i$ are positive integers, we always have
$$
n_1^2+n_2^2+...n_k^2<n^2.
$$ 

This means in the summation, the term of $k=1$ will be the leading term and will determine its order. Thus we have
$$
\wh\b_{X/\F_q,n}(0)=O(\whv_n)=O(q^{(g-1)n(n+1)/2})
$$
as wanted. The rest parts are trivial.
\ep

\begin{cor}
In particular, when $n\gg0$,
$$
\frac{\wh\b_{X/\F_q,n}(0)}{\whv_n}=1+o(1).
$$
\end{cor}

\section{High Rank Zeta Function over Genus 2 Curves}
\subsection{Multiplicative Structure of Artin Zeta Function}
As a preparation, we will recall the multiplicative structure of Artin zeta function. This is equivalent to the Euler product structure (\ref{EPF}) for Artin zeta function. Namely, by applying the rationality and functional equation, we have 
$$\ba
\zeta_{X/\F_q}(s)=&\frac{\prod_{i=1}^g(1-\om_{X/\F_q,i}t)(1-\om^{'}_{X/\F_q,i}t)}{(1-t)(1-qt)}
\\=&\exp\Bigg(\sum_{i=1}^g\log(1-\om_{X/\F_q,i}t)+\log(1-\om^{'}_{X/\F_q,i}t)
\\&\ \ \ \ \ \ \ \ -\log(1-t)-\log(1-qt)\Bigg)
\\=&\exp\left(\sum_{k=1}^\infty N_k\frac{t^k}{k}\right),
\ea
$$
where 
$$
N_k=q^k+1-\sum_{i=1}^g\left(\om_{X/\F_q,i}^k+{\om'}_{X/\F_q,i}^k\right).
$$
One may easily find that 
$$
N_1=\#X(\F_q),
$$
and furthermore, according to (\ref{EPF}), we could conclude that 
$$
N_k=\#X(\F_{q^k}).
$$

In fact, this is exactly another definition of zeta function for curves over finite field, which is called Hasse--Weil zeta function, generally defined for algebraic varieties over finite fields. Based on the above discussion, Artin zeta function and Hasse--Weil zeta function coincide with each other for $X/\F_q$.
\begin{mdef}[Hasse--Weil Zeta Function]
For a(n integral regular projective) curve $X/\F_q$ defined over the finite field $\F_q$ with $q$ elements, its Hasse--Weil zeta function is defined by
$$
\zeta_{X/\F_q}(s):=\exp\left(\sum_{k=1}^\infty\frac{\#X(\F_{q^k})}{k}(q^{-s})^k\right)\qquad\Re(s)>1,
$$
where $X(\F_{q^k})$ denotes the set of $\F_{q^k}$-rational points of $X$. 
\end{mdef}

\subsection{$P_{X/\F_q,n}(T)$ for Genus 2 Curves}
In this subsection, we will fix $X/\F_q$ as a genus $g=2$ curve over finite field $\F_q$, unless there is an additional assumption. We start with $n=1$, then rank $n$ zeta function $\zeta_{X/\F_q,n}(s)$ for genus 2 curve $X/\F_q$ will reduce to Artin zeta function, or equivalently, Hasse--Weil zeta function $\zeta_{X/\F_q}(s)$. In particular, based on Theorem\ \ref{AZP}(1), we have 
$$
\ba
\zeta_{X/\F_q}(s)=&\frac{P_{X/\F_q}(t)}{(1-t)(1-qt)}
\\=&\frac{1+a_{X/\F_q}t+b_{X/\F_q}t^2+a_{X/\F_q}qt^3+q^2t^4}{(1-t)(1-qt)}
\\=&\exp\left(\sum_{k=1}^\infty\frac{\#X(\F_{q^k})}{k}(q^{-s})^k\right).
\ea
$$ 
Then after taking logarithms on both sides, differentiating with respect to $T$ and comparing the coefficients, we may conclude that the $N_1$ and $N_2$ suffice to determine the coefficients $a_{X/\F_q}$ and $b_{X/\F_q}$, and the result is that:
$$
\bc
a_{X/\F_q}=N_1-1-q,
\\
b_{X/\F_q}=\frac{1}{2}\cdot\left(N_2-1-q^2+a_{X/\F_q}^2\right),
\ec
$$
where $N_{i}=\#X(\F_{q^i})$ denotes the number of $\F_{q^i}$-rational points of $X$ for $i\in\{1,\ 2\}$, as just introduced. This means the coefficients $a_{X/\F_q}$ and $b_{X/\F_q}$ have a strong geometric meaning since they have a strong relation to number of rational points of the curve $X$.

Next, we may focus on $P_{X/\F_q,n}(s)$ for general $n\ge 1$. Recall that in the equation \ref{PT}, we have known that
$$
\ba
&P_{X/\F_q,n}(T)=\wh Z_{X/\F_q,n}(T)\cdot(1-T)(1-QT)\cdot T^{g-1}
\\=&\Big(\sum_{m=0}^{g-2}\a_{X/\F_q,n}(nm)(T^m+Q^{g-1-m}T^{2g-2-m})+\a_{X/\F_q,n}(n(g-1))T^{g-1}\Big)
\\&\times(1-T)(1-QT)+(Q-1)T^g\b_{X/\F_q,n}(0).
\ea
$$
Then suppose $g=2$, we have:
$$
\begin{aligned}
P_{X/\F_q,n}(T)=
&\left(\a_{X/\F_q,n}(0)(1+QT^2)+\a_{X/\F_q,n}(n)T\right)
\\&\times(1-T)(1-QT)+(Q-1)T^2\b_{X/\F_q,n}(0)
\\=&\Bigg(1+\left(\frac{\a_{X/\F_q,n}(n)}{\a_{X/\F_q,n}(0)}-(Q+1)\right)T
\\&+\left((Q-1)\frac{\b_{X/\F_q,n}(0)}{\a_{X/\F_q,n}(0)}+2Q-(Q+1)\frac{\a_{X/\F_q,n}(n)}{\a_{X/\F_q,n}(0)}\right)T^2
\\&+Q\left(\frac{\a_{X/\F_q,n}(n)}{\a_{X/\F_q,n}(0)}-(Q+1)\right)T^3
\\&+Q^2T^4\Bigg)\a_{X/\F_q,n}(0).
\end{aligned}
$$
Accordingly, we may obtain
$$
\frac{P_{X/\F_q,n}(T)}{\a_{X/\F_q,n}(0)}=1+a_{X/\F_q,n}T+b_{X/\F_q,n}T^2+a_{X/\F_q,n}QT^3+Q^2T^4,
$$
where we have set
\be\label{anbn}
\bc
a_{X/\F_q,n}:=\frac{\a_{X/\F_q,n}(n)}{\a_{X/\F_q,n}(0)}-(Q+1),
\\
b_{X/\F_q,n}:=(Q-1)\frac{\b_{X/\F_q,n}(0)}{\a_{X/\F_q,n}(0)}+2Q-(Q+1)\frac{\a_{X/\F_q,n}(n)}{\a_{X/\F_q,n}(0)}.
\ec
\ee
One may easily check that, based on the naturality, that is, Theorem\ \ref{HRZP}(1), we have
$$
\bc
a_{X/\F_q,1}=a_{X/\F_q}=\frac{\a_{X/\F_q,1}(1)}{\a_{X/\F_q,1}(0)}-(q+1),
\\
b_{X/\F_q,1}=b_{X/\F_q}=(q-1)\frac{\b_{X/\F_q,1}(0)}{\a_{X/\F_q,1}(0)}+2q-(q+1)\frac{\a_{X/\F_q,1}(1)}{\a_{X/\F_q,1}(0)}.
\ec
$$
Since by definition, 
$$
\a_{X/\F_q,1}(0)=\sum_{\cL\in\Pic^0(X)}\frac{q^{h^0(X,\cL)}-1}{q-1}=\frac{q^{h^0(X,\cO_X)}-1}{q-1}=1,
$$
we could obtain the following from the first equation:
\be\label{alpha_1_1}
\a_{X/\F_q,1}(1)=a_{X/\F_q,1}+q+1=N_1.
\ee
This will be used later in our computation. One may find that the definition of $\a_{X/\F_q,1}(1)$ and definition of $N_1$ are equivalent. On the other hand, the second equation is, in fact, equivalent to the property of residue in geometry, Theorem\ \ref{HRZP}(4), that is,
$$
\beta_{X/\F_q,1}(0)=\Res_{t=1}Z_{X/\F_q,1}(t)=\frac{1+a_{X/\F_q}+b_{X/\F_q}+a_{X/\F_q}q+q^2}{q-1}.
$$

We end this section with the following comments. Based on several geometric theorems and zeta properties, rank $n$ zeta function $\zeta_{X/\F_q,n}(s)$ over a genus $2$ curve $X/\F_q$ can be expressed into a relatively simple form, involving two rank $n$ coefficients defined in (\ref{anbn}), $a_{X/\F_q,n}$ and $b_{X/\F_q,n}$, while these two rational numbers are totally determined by the geo-arithmetic invariants $\a_{X/\F_q,n}(n)$, $\a_{X/\F_q,n}(0)$ and $\b_{X/\F_q,n}(0)$. By examining the properties of $\a$ and $\b$ invariants carefully, we will discover some structures of the rank $n$ zeta function $\zeta_{X/\F_q,n}$, so that we could asymptotically prove the Riemann hypothesis of a stronger form.

\section{Riemann Hypothesis for $\wh\zeta_{X/\F_q,n}$ over Genus 2 Curves}
In this section, we establish the asymptotic Riemann Hypothesis for the case $g=2$. The main proposition depends on a detailed examination of the structural properties of the associated invariants $\a_{X/\F_q,n}(n)$, $\a_{X/\F_q,n}(0)$ and $\b_{X/\F_q,n}(0)$. The initial motivation arose from computational investigations: numerical experiments for several examples of genus g = 2 curves indicated that the coefficients of the corresponding zeta functions tend to approach $\pm\frac{\sqrt{2}}{2}$. These computational results suggested that, for $g=2$, the Riemann Hypothesis holds asymptotically as $n\to\infty$, and moreover, that the distribution of zeros exhibits a structure even subtler than that predicted by the Riemann hypothesis.

Recall that in the previous section, we have verified in equation (\ref{anbn}) that
$$
\frac{P_{X/\F_q,n}(T)}{\a_{X/\F_q,n}(0)}=1+a_{X/\F_q,n}T+b_{X/\F_q,n}T^2+a_{X/\F_q,n}QT^3+Q^2T^4,
$$
where we have set
$$
\bc
a_{X/\F_q,n}:=\frac{\a_{X/\F_q,n}(n)}{\a_{X/\F_q,n}(0)}-(Q+1),
\\
b_{X/\F_q,n}:=(Q-1)\frac{\b_{X/\F_q,n}(0)}{\a_{X/\F_q,n}(0)}+2Q-(Q+1)\frac{\a_{X/\F_q,n}(n)}{\a_{X/\F_q,n}(0)}.
\ec
$$
We will prove that these two coefficients converge to constant numbers as $n\to\infty$, which will lead to a better bound than Riemann Hypothesis for large rank case.

A central difficulty encountered in the development of this idea was the analysis of the $\a_{X/\F_q,n}(n)$ invariant. Motivated by Weng, this difficulty was resolved by referring to \cite{RH}, where the \textit{General Counting Miracle Theorem} 
was introduced. That theorem provides a closed formula for $\a_{X/\F_q,n}(n)$, which is a key component of the present argument. This was derived from Theorem\,\ref{RoHRZ} and Theorem\,\ref{SU}. In what follows, we will focus on the case $g = 2$; however, the reference \cite{RH} treats the general case for arbitrary $g$. For additional details, one may refer to \cite{RH}.

\begin{prop}[A Closed Formula for $\a_{X/\F_q;n}(n)$, see \cite{RH} for general curve]
Assume that $g=2$. Then we can obtain 
$$
\ba
&q^{-\binom{n}{2}}\a_{X/\F_q,n}(n)
\\&=\sum_{a=1}^n\Res_{T=0}T^{-1}\Biggl(\sum_{\substack{k_1,\ldots,k_p>0\\ k_1+\ldots+k_p=n-a}}\frac{\wh v_{k_1}\ldots\wh v_{k_p}}{\prod_{j=1}^{p-1}(1-q^{k_j+k_{j+1}})} (-1)\sum_{\ell=1}^\infty\bl q^{n-a-k_{p}}T\br^\ell
\\&\times\left(\a_{X/\F_q,1}(0)\left(q^{-(n-a)}T^{-1}\right)+\a_{X/\F_q,1}\bl1\br\right)
\\&\times\sum_{\substack{l_1,\ldots,l_r>0\\ l_1+\ldots+l_r=a-1}}
 \frac{\wh v_{l_1}\ldots\wh v_{l_r}}{\prod_{j=1}^{r-1}(1-q^{l_j+l_{j+1}})}\sum_{\kappa=0}^\infty\bl q^{n-a+1+l_{1}}T\br^\kappa\Biggr).
\ea
$$
\end{prop}

This means we must consider the constant term inside the bracket.
For the term consists $\a_{X/\F_q,1}(1)$,  since the first partition starts from $T$, the second partition starts from $1$, it is non-zero if and only if the first partition is $1$, that is, $a=n$. 
That is, 

$$
\ba
&\sum_{a=1}^n\Res_{T=0}T^{-1}\Biggl(\sum_{\substack{k_1,\ldots,k_p>0\\ k_1+\ldots+k_p=n-a}}\frac{\wh v_{k_1}\ldots\wh v_{k_p}}{\prod_{j=1}^{p-1}(1-q^{k_j+k_{j+1}})} (-1)\sum_{\ell=1}^\infty\bl q^{n-a-k_{p}}T\br^\ell
\\&\times\left(\a_{X/\F_q,1}\bl1\br\right)
\\&\times\sum_{\substack{l_1,\ldots,l_r>0\\ l_1+\ldots+l_r=a-1}}
\frac{\wh v_{l_1}\ldots\wh v_{l_r}}{\prod_{j=1}^{r-1}(1-q^{l_j+l_{j+1}})}\sum_{\kappa=0}^\infty\bl q^{n-a+1+l_{1}}T\br^\kappa\Biggr)
\\&=\Res_{T=0}T^{-1}\Biggl(1\times\left(\a_{X/\F_q,1}\bl1\br\right)\times\sum_{\substack{l_1,\ldots,l_r>0\\ l_1+\ldots+l_r=n-1}}
 \frac{\wh v_{l_1}\ldots\wh v_{l_r}}{\prod_{j=1}^{r-1}(1-q^{l_j+l_{j+1}})}\sum_{\kappa=0}^\infty\bl q^{1+l_{1}}T\br^\kappa\Biggr)
\\&=\a_{X/\F_q,1}(1)\times\sum_{\substack{l_1,\ldots,l_r>0\\ l_1+\ldots+l_r=n-1}}
 \frac{\wh v_{l_1}\ldots\wh v_{l_r}}{\prod_{j=1}^{r-1}(1-q^{l_j+l_{j+1}})}
\\&=\a_{X/\F_q,1}(1)\times\wh\beta_{X/\F_q,n-1}(0).
\ea
$$

As for the term $\a_{X/\F_q,1}(0)\left(q^{-(n-a)}T^{-1}\right)$, the summation can be decomposed into three part:
$$
\sum_{a=1}^n=\sum_{a=n}+\sum_{a=1}+\sum_{a=2}^{n-1}.
$$
For the part $a=n$, the first partition is 1.
Thus, it is
$$
\sum_{a=n}\Res_{T=0}T^{-1}\Biggl(1\times\left(T^{-1}\right)\times\sum_{\substack{l_1,\ldots,l_r>0\\ l_1+\ldots+l_r=n-1}}
 \frac{\wh v_{l_1}\ldots\wh v_{l_r}}{\prod_{j=1}^{r-1}(1-q^{l_j+l_{j+1}})}\sum_{\kappa=0}^\infty\bl q^{1+l_{1}}T\br^\kappa\Biggr).
$$
The constant term is the case when $\kappa=1$.
That is,
$$
\sum_{\substack{l_1,\ldots,l_r>0\\ l_1+\ldots+l_r=n-1}}
 \frac{\wh v_{l_1}\ldots\wh v_{l_r}}{\prod_{j=1}^{r-1}(1-q^{l_j+l_{j+1}})}\bl q^{1+l_{1}}\br.
$$
Similarly, for the part $a=1$, the second partition is 1. That is,
$$
\ba
\sum_{a=1}\Res_{T=0}T^{-1}&\Biggl(\sum_{\substack{k_1,\ldots,k_p>0\\ k_1+\ldots+k_p=n-1}}\frac{\wh v_{k_1}\ldots\wh v_{k_p}}{\prod_{j=1}^{p-1}(1-q^{k_j+k_{j+1}})} (-1)\sum_{\ell=1}^\infty\bl q^{n-1-k_{p}}T\br^\ell
\\&\times\left(q^{-(n-1)}T^{-1}\right)\times1\Biggr).
\ea
$$
The constant term is the case when $\ell=1$. That is,
$$
(-1)\sum_{\substack{k_1,\ldots,k_p>0\\ k_1+\ldots+k_p=n-1}}\frac{\wh v_{k_1}\ldots\wh v_{k_p}}{\prod_{j=1}^{p-1}(1-q^{k_j+k_{j+1}})}\bl q^{-k_{p}}\br.
$$
For the part $\sum_{a=2}^{n-1}$, these two ordered partitions do not degenerate:
$$
\ba
&\sum_{a=2}^{n-1}\Res_{T=0}T^{-1}\Biggl(\sum_{\substack{k_1,\ldots,k_p>0\\ k_1+\ldots+k_p=n-a}}\frac{\wh v_{k_1}\ldots\wh v_{k_p}}{\prod_{j=1}^{p-1}(1-q^{k_j+k_{j+1}})} (-1)\sum_{\ell=1}^\infty\bl q^{n-a-k_{p}}T\br^\ell
\\&\times\left(\left(q^{-(n-a)}T^{-1}\right)\right)
\\&\times\sum_{\substack{l_1,\ldots,l_r>0\\ l_1+\ldots+l_r=a-1}}
 \frac{\wh v_{l_1}\ldots\wh v_{l_r}}{\prod_{j=1}^{r-1}(1-q^{l_j+l_{j+1}})}\sum_{\kappa=0}^\infty\bl q^{n-a+1+l_{1}}T\br^\kappa\Biggr).
\ea
$$
To obtain constant term, $\ell$ must be $1$ and $\kappa$ must be $0$. That is,
$$
\ba
&\sum_{a=2}^{n-1}\Biggl(\sum_{\substack{k_1,\ldots,k_p>0\\ k_1+\ldots+k_p=n-a}}\frac{\wh v_{k_1}\ldots\wh v_{k_p}}{\prod_{j=1}^{p-1}(1-q^{k_j+k_{j+1}})} (-1)\bl q^{n-a-k_{p}}\br
\\&\times\left(q^{-(n-a)}\right)\times\sum_{\substack{l_1,\ldots,l_r>0\\ l_1+\ldots+l_r=a-1}}
 \frac{\wh v_{l_1}\ldots\wh v_{l_r}}{\prod_{j=1}^{r-1}(1-q^{l_j+l_{j+1}})}\Biggr)
\\&=\sum_{a=2}^{n-1}\Biggl(\sum_{\substack{k_1,\ldots,k_p>0\\ k_1+\ldots+k_p=n-a}}\frac{\wh v_{k_1}\ldots\wh v_{k_p}}{\prod_{j=1}^{p-1}(1-q^{k_j+k_{j+1}})} (-1)\bl q^{-k_{p}}\br
\\&\times\sum_{\substack{l_1,\ldots,l_r>0\\ l_1+\ldots+l_r=a-1}}
 \frac{\wh v_{l_1}\ldots\wh v_{l_r}}{\prod_{j=1}^{r-1}(1-q^{l_j+l_{j+1}})}\Biggr).
\ea
$$

Finally, we summarize that 
\be\label{alphann}
\ba
q^{-\binom{n}{2}}\a_{X/\F_q,n}(n)&=\a_{X/\F_q,1}(1)\times\wh\beta_{X/\F_q,n-1}(0)
\\&+\sum_{\substack{l_1,\ldots,l_r>0\\ l_1+\ldots+l_r=n-1}}
\frac{\wh v_{l_1}\ldots\wh v_{l_r}}{\prod_{j=1}^{r-1}(1-q^{l_j+l_{j+1}})}\bl q^{1+l_{1}}\br
\\&+(-1)\sum_{\substack{k_1,\ldots,k_p>0\\ k_1+\ldots+k_p=n-1}}\frac{\wh v_{k_1}\ldots\wh v_{k_p}}{\prod_{j=1}^{p-1}(1-q^{k_j+k_{j+1}})}\bl q^{-k_{p}}\br
\\&+\sum_{a=2}^{n-1}\Biggl(\sum_{\substack{k_1,\ldots,k_p>0\\ k_1+\ldots+k_p=n-a}}\frac{\wh v_{k_1}\ldots\wh v_{k_p}}{\prod_{j=1}^{p-1}(1-q^{k_j+k_{j+1}})} (-1)\bl q^{-k_{p}}\br
\\&\times\sum_{\substack{l_1,\ldots,l_r>0\\ l_1+\ldots+l_r=a-1}}
 \frac{\wh v_{l_1}\ldots\wh v_{l_r}}{\prod_{j=1}^{r-1}(1-q^{l_j+l_{j+1}})}\Biggr).
\ea
\ee

Next, we will introduce the main proposition of this section. It gives the estimate of $\frac{\b_{X/\F_q,n}(0)}{\a_{X/\F_q,n}(0)}$ and $\frac{\a_{X/\F_q,n}(n)}{\a_{X/\F_q,n}(0)}$, which appear in the polynomial $\frac{P_{X/\F_q,n}(T)}{\a_{X/\F_q,n}(0)}$.

\begin{prop}[Estimate of $\frac{\b_{X/\F_q,n}(0)}{\a_{X/\F_q,n}(0)}$ and $\frac{\a_{X/\F_q,n}(n)}{\a_{X/\F_q,n}(0)}$]\label{MP}
For a genus $2$ curve $X/\F_q$, we have, when $n\to\infty$,
\begin{itemize}
\item[(1)] $\frac{\b_{X/\F_q,n}(0)}{\a_{X/\F_q,n}(0)}\to Q+N_1;$
\item[(2)] $\frac{\a_{X/\F_q,n}(n)}{\a_{X/\F_q,n}(0)}\to Q+N_1;$
\item[(3)] $Q\frac{\b_{X/\F_q,n}(0)}{\a_{X/\F_q,n}(0)}\to Q^2+N_1Q+q+(q-1)\whv_1+2q^2\whv_1-2\whv_1;$
\item[(4)] $Q\frac{\a_{X/\F_q,n}(n)}{\a_{X/\F_q,n}(0)}\to Q^2+N_1Q+2q^2\whv_1-q\whv_1-\whv_1-q$.
\end{itemize}
where $N_i:=\#X(\F_{q^i})$ denotes the number of $\F_{q^i}$-rational points of $X$.
\end{prop}

To prove this proposition, we first mention that (1) and (2) are direct conclusions from (3) and (4) respectively, by dividing $Q$ on both sides. This means we only have to prove (3) and (4).

Firstly, for $Q\frac{\b_{X/\F_q,n}(0)}{\a_{X/\F_q,n}(0)}$, according to definition and Theorem\,\ref{CM}, we have
$$
\ba
Q\frac{\beta_{X/\F_q,n}(0)}{\a_{X/\F_q,n}(0)}&=Q\frac{q^{n(n-1)/2}\wh\b_{X/\F_q,n}(0)}{q^{n-1}\b_{X/\F_q,n-1}(0)}
\\&=Q\frac{q^{n(n-1)/2}\wh\b_{X/\F_q,n}(0)}{q^{n-1}q^{(n-1)(n-2)/2}\wh\b_{X/\F_q,n-1}(0)}
\\&=Q\frac{\wh\b_{X/\F_q,n}(0)}{\wh\b_{X/\F_q,n-1}(0)}.
\ea
$$

Thus, according to Theorem\,\ref{SSM}, we obtain
$$
\ba
Q\frac{\b_{X/\F_q,n}(0)}{\a_{X/\F_q,n}(0)}&=Q\frac{\wh\b_{X/\F_q,n}(0)}{\wh\b_{X/\F_q,n-1}(0)}
\\&=Q\frac{\whv_{n}}{\wh\b_{X/\F_q,n-1}(0)}-Q\sum_{i=1}^{n-1}\frac{\whv_i\whv_{n-i}}{(Q-1)\wh\b_{X/\F_q,n-1}(0)}+o(1).
\ea
$$
In this equation, we first apply Theorem\,\ref{SSM} to $\wh\b_{X/\F_q,n}(0)$. The summation is divided into three parts, based on the length of the ordered partition on $n$. It is easy to check that if the length of the partition of $n$ is lager than 2, the summation is o(1), which means it goes to 0 when n goes to infinity.

For the first part, we apply Theorem\,\ref{SSM} to $\wh\b_{X/\F_q,n-1}(0)$ again to obtain
$$
\ba
Q\frac{\whv_{n}}{\wh\b_{X/\F_q,n-1}(0)}&=Q\frac{\whv_{n}}{\whv_{n-1}-\sum_{i=1}^{n-2}\frac{\whv_{i}\whv_{n-1-i}}{q^{n-1}-1}+...}
\\&=Q\frac{\whv_{n}}{\whv_{n-1}}\times\frac{1}{1-\sum_{i=1}^{n-2}\frac{\whv_{i}\whv_{n-1-i}}{(q^{n-1}-1)\whv_{n-1}}+...}
\\&\to Q\wh\z_{X/\F_q}(n)\times\frac{1}{1-2\frac{\whv_{1}\whv_{n-2}}{(q^{n-1}-1)\whv_{n-1}}}
\\&= Q\wh\z_{X/\F_q}(n)+2Q\frac{\whv_{1}\wh\z_{X/\F_q}(n)}{(q^{n-1}-1)\wh\z_{X/\F_q}(n-1)}+o(1)
\\&\to Q\wh\z_{X/\F_q}(n)+2q\whv_{1}\times\frac{\wh\z_{X/\F_q}(n)}{\wh\z_{X/\F_q}(n-1)}\times\frac{1}{1-q^{-n+1}}.
\ea
$$

For the second part, we have 
$$
-Q\sum_{i=1}^{n-1}\frac{\whv_i\whv_{n-i}}{(Q-1)\wh\b_{X/\F_q,n-1}(0)}=-\sum_{i=1}^{n-1}\frac{\whv_i\whv_{n-i}}{\wh\b_{X/\F_q,n-1}(0)}\times\frac{1}{1-q^n} \to-2\whv_1.
$$
In the above estimations, we have use the fact that $\frac{1}{1-x}=1+x+x^2+...$, together with the estimate of $\whv_n$ and $\wh\b_{X/\F_q,n}(0)$, introduced in section 3. 

In the previous version of this study, we use big O-notation to estimate these coefficients. This means we will only consider terms of order $O(q^{>0})$ and ignore all finite constant terms. Even if we could obtain the same result, it should be more suitable to use small o-notation do the estimation, so that those constant terms will be remained.

Thus, we introduce the following lemma:
\begin{lem}
For a fixed genus $2$ curve $X/\F_q$, we have:
\begin{itemize}
\item[(1)] $Q\wh\z_{X/\F_q}(n)\to Q^2+N_1Q+q+(q-1)\whv_1;$
\item[(2)] $\frac{\wh\z_{X/\F_q}(n)}{\wh\z_{X/\F_q}(n-1)}\to q$.
\end{itemize}
\end{lem}
\bp
These estimations come from the definition of $\wh\z_{X/\F_q}(n)$ and the the geometric series expansion of $\frac{1}{1-x}$. The main reason for these constant terms is the shift $Q=q^n$.

Recall that for a genus $2$ curve $X/\F_q$,
$$
\whv_1=\wh\z_{X/\F_q}^*(1)=\frac{q^2+a_{X/\F_q}q+a_{X/\F_q}+b_{X/\F_q}+1}{q-1},
$$
$$
\wh\z_{X/\F_q}(n)=\frac{q^n+a_{X/\F_q}+\frac{b_{X/\F_q}}{q^n}+\frac{a_{X/\F_q}}{q^{2n-1}}+\frac{1}{q^{3n-2}}}{(1-q^{-n})(1-q^{1-n})}.
$$
Thus, (2) is trivial.

As for (1), we have the following estimation
$$
\ba
Q\wh\z_{X/\F_q}(n)=&Q\frac{q^n+a_{X/\F_q}+\frac{b_{X/\F_q}}{q^n}+\frac{a_{X/\F_q}}{q^{2n-1}}+\frac{1}{q^{3n-2}}}{(1-q^{-n})(1-q^{1-n})}
\\=&Q\bl q^n+a_{X/\F_q}+\frac{b_{X/\F_q}}{q^n}+\frac{a_{X/\F_q}}{q^{2n-1}}+\frac{1}{q^{3n-2}}\br
\\&\times(1+q^{-n}+q^{-2n}+q^{-3n}+...)
\\&\times(1+q^{1-n}+q^{2-2n}+q^{3-3n}+...)
\\=&\bl q^n+a_{X/\F_q}+\frac{b_{X/\F_q}}{q^n}+\frac{a_{X/\F_q}}{q^{2n-1}}+\frac{1}{q^{3n-2}}\br
\\&\times(q^n+1+q^{-n}+q^{-2n}+...)
\\&\times(1+q^{1-n}+q^{2-2n}+q^{3-3n}+...)
\\ \to& q^{2n}+q^{n+1}+q^{2}+q^n+q+1+q^na_{X/\F_q}+qa_{X/\F_q}+a_{X/\F_q}+b_{X/\F_q}
\\=&q^{2n}+q^n(q+1+a_{X/\F_q})+(q^2+a_{X/\F_q}q+a_{X/\F_q}+b_{X/\F_q}+1)+q
\\=&Q^2+N_1Q+(q-1)\whv_1+q
\ea
$$
as wanted.
\ep

Based on the above lemma, we could summarize that 
$$
\ba
Q\frac{\b_{X/\F_q,n}(0)}{\a_{X/\F_q,n}(0)}&\to Q\frac{\whv_{n}}{\wh\b_{X/\F_q,n-1}(0)}-Q\sum_{i=1}^{n-1}\frac{\whv_i\whv_{n-i}}{(Q-1)\wh\b_{X/\F_q,n-1}(0)}
\\&\to Q\wh\z_{X/\F_q}(n)+2q\whv_{1}\times\frac{\wh\z_{X/\F_q}(n)}{\wh\z_{X/\F_q}(n-1)}\times\frac{1}{1-q^{-n+1}}-2\whv_1
\\&\to Q^2+N_1Q+(q-1)\whv_1+q+2q^2\whv_1-2\whv_1.
\ea
$$
This proves Proposition\,\ref{MP} (3), and hence, Proposition\,\ref{MP} (1).

On the other hand, for $Q\frac{\a_{X/\F_q,n}(n)}{\a_{X/\F_q,n}(0)}$, based on the definition and Theorem\,\ref{CM}, we first summarize the relation:
$$
\ba
\frac{\a_{X/\F_q,n}(n)}{\a_{X/\F_q,n}(0)}&=\frac{\a_{X/\F_q,n}(n)}{q^{n-1}\b_{X/\F_q,n-1}(0)}
\\&=\frac{\a_{X/\F_q,n}(n)}{q^{n-1}q^{(n-1)(n-2)/2}\wh\b_{X/\F_q,n-1}(0)}
\\&=\frac{q^{-n(n-1)/2}\a_{X/\F_q,n}(n)}{\wh\b_{X/\F_q,n-1}(0)}.
\ea
$$

Thus, according to equation (\ref{alphann}), we have 
$$
\ba
Q\frac{\a_{X/\F_q,n}(n)}{\a_{X/\F_q,n}(0)}&=Q\Biggl(\a_{X/\F_q,1}(1)\times\wh\beta_{X/\F_q,n-1}(0)
\\&+\sum_{\substack{l_1,\ldots,l_r>0\\ l_1+\ldots+l_r=n-1}}
\frac{\wh v_{l_1}\ldots\wh v_{l_r}}{\prod_{j=1}^{r-1}(1-q^{l_j+l_{j+1}})}\bl q^{1+l_{1}}\br
\\&+(-1)\sum_{\substack{k_1,\ldots,k_p>0\\ k_1+\ldots+k_p=n-1}}\frac{\wh v_{k_1}\ldots\wh v_{k_p}}{\prod_{j=1}^{p-1}(1-q^{k_j+k_{j+1}})}\bl q^{-k_{p}}\br
\\&+\sum_{a=2}^{n-1}\Biggl(\sum_{\substack{k_1,\ldots,k_p>0\\ k_1+\ldots+k_p=n-a}}\frac{\wh v_{k_1}\ldots\wh v_{k_p}}{\prod_{j=1}^{p-1}(1-q^{k_j+k_{j+1}})} (-1)\bl q^{-k_{p}}\br
\\&\times\sum_{\substack{l_1,\ldots,l_r>0\\ l_1+\ldots+l_r=a-1}}
 \frac{\wh v_{l_1}\ldots\wh v_{l_r}}{\prod_{j=1}^{r-1}(1-q^{l_j+l_{j+1}})}\Biggr)\Biggr)/\wh\b_{X/\F_q,n-1}(0).
\ea
$$

Then, we will consider the asymptotic behavior of each term in the equation.
\begin{itemize}
\item[(1)] $Q\bl\a_{X/\F_q,1}(1)\times\wh\beta_{X/\F_q,n-1}(0)\br/\wh\beta_{X/\F_q,n-1}(0)$.

Easily, this is just $Q\a_{X/\F_q,1}(1)=QN_1$.

\item[(2)] $Q\bl\sum_{\substack{l_1,\ldots,l_r>0\\ l_1+\ldots+l_r=n-1}}
\frac{\wh v_{l_1}\ldots\wh v_{l_r}}{\prod_{j=1}^{r-1}(1-q^{l_j+l_{j+1}})}\bl q^{1+l_{1}}\br\br/\wh\beta_{X/\F_q,n-1}(0)$.

We rewrite the summation of ordered partition by the length $r$:
$$
\ba
&Q\bl\sum_{\substack{l_1,\ldots,l_r>0\\ l_1+\ldots+l_r=n-1}}
\frac{\wh v_{l_1}\ldots\wh v_{l_r}}{\prod_{j=1}^{r-1}(1-q^{l_j+l_{j+1}})}\bl q^{1+l_{1}}\br\br/\wh\beta_{X/\F_q,n-1}(0)
\\&=Q\frac{q^n\whv_{n-1}}{\wh\b_{X/\F_q,n-1}(0)}-Q\sum_{i=1}^{n-2}\frac{q^{1+i}\whv_i\whv_{n-1-i}}{(q^{n-1}-1)\wh\b_{X/\F_q,n-1}(0)}+o(1)
\\&\to Q\frac{q^n\whv_{n-1}}{\wh\b_{X/\F_q,n-1}(0)}-\frac{q^{n+1}}{q^{n-1}-1}\sum_{i=1}^{n-2}\frac{q^i\whv_i\whv_{n-1-i}}{\wh\b_{X/\F_q,n-1}(0)}
\\&\to Q^2\frac{1}{1-2\frac{\whv_1\whv_{n-2}}{(q^{n-1}-1)\whv_{n-1}}}-(q^2+q^{3-n}+q^{4-2n}+...)\sum_{i=1}^{n-2}\frac{q^i\whv_i\whv_{n-1-i}}{\wh\b_{X/\F_q,n-1}(0)}
\\&\to Q^2(1+2\frac{\whv_1\whv_{n-2}}{(q^{n-1}-1)\whv_{n-1}}+...)-q^2\frac{q^{n-2}\whv_1\whv_{n-2}}{\wh\b_{X/\F_q,n-1}(0)}
\\&\to Q^2+2\whv_1\frac{Q^2}{(q^{n-1}-1)\wh\z_{X/\F_q}(n-1)}-Q\frac{\whv_1\whv_{n-2}}{\whv_{n-1}-\sum_{i=1}^{n-2}\frac{\whv_i\whv_{n-i-1}}{q^{n-1}-1}+...}
\\&\to Q^2+2q^2\whv_1-Q\frac{\whv_1}{\wh\z_{X/\F_q}(n-1)}
\\&\to Q^2+2q^2\whv_1-q\whv_1.
\ea
$$

\item[(3)] $(-1)Q\Biggl(\sum_{\substack{k_1,\ldots,k_p>0\\ k_1+\ldots+k_p=n-1}}\frac{\wh v_{k_1}\ldots\wh v_{k_p}}{\prod_{j=1}^{p-1}(1-q^{k_j+k_{j+1}})}\bl q^{-k_{p}}\br\Biggr)/\wh\b_{X/\F_q,n-1}(0)$.

Based on the same reason, we have the following estimation:
$$
\ba
&-Q\Biggl(\sum_{\substack{k_1,\ldots,k_p>0\\ k_1+\ldots+k_p=n-1}}\frac{\wh v_{k_1}\ldots\wh v_{k_p}}{\prod_{j=1}^{p-1}(1-q^{k_j+k_{j+1}})}\bl q^{-k_{p}}\br\Biggr)/\wh\b_{X/\F_q,n-1}(0)
\\&=-Q\frac{\whv_{n-1}}{q^{n-1}}/\wh\b_{X/\F_q,n-1}(0)+Q\sum_{i=1}^{n-2}\frac{\whv_i\whv_{n-i-1}}{(q^{n-1}-1)q^i\wh\b_{X/\F_q,n-1}(0)}+...
\\&\to-\frac{Q}{q^{n-1}}\times\frac{\whv_{n-1}}{\wh\b_{X/\F_q,n-1}(0)}
\\&\to-q.
\ea
$$

\item[(4)] For the rest part, we have
$$
\ba
&Q\sum_{a=2}^{n-1}\Biggl(\sum_{\substack{k_1,\ldots,k_p>0\\ k_1+\ldots+k_p=n-a}}\frac{\wh v_{k_1}\ldots\wh v_{k_p}}{\prod_{j=1}^{p-1}(1-q^{k_j+k_{j+1}})} (-1)\bl q^{-k_{p}}\br
\\ &\times\sum_{\substack{l_1,\ldots,l_r>0\\ l_1+\ldots+l_r=a-1}}
 \frac{\wh v_{l_1}\ldots\wh v_{l_r}}{\prod_{j=1}^{r-1}(1-q^{l_j+l_{j+1}})}\Biggr)/\wh\b_{X/\F_q,n-1}(0)
\\&=-Q\sum_{a=2}^{n-1}\Biggl(\sum_{\substack{k_1,\ldots,k_p>0\\ k_1+\ldots+k_p=n-a}}\frac{\wh v_{k_1}\ldots\wh v_{k_p}}{\prod_{j=1}^{p-1}(1-q^{k_j+k_{j+1}})}\bl q^{-k_{p}}\br\times\wh\b_{X/\F_q,a-1}(0)\Biggr)/\wh\b_{X/\F_q,n-1}(0)
\\&\to-\whv_1.
\ea
$$
We mention that for $a\ne n-1$, the summation goes to $0$; and when $a=n-1$, it goes to $-\whv_1$ when $n\to\infty$. 

\end{itemize}

Finally, we could summarize that 
$$
Q\frac{\a_{X/\F_q,n}(n)}{\a_{X/\F_q,n}(0)}\to QN_1+Q^2+2q^2\whv_1-q\whv_1-q-\whv_1
$$
as wanted. This proves Proposition\,\ref{MP} (4), and hence, Proposition\,\ref{MP} (2).

\begin{prop}[Estimate of $a_{X/\F_q,n}$ and $b_{X/\F_q,n}$]
For a genus $2$ curve $X/\F_q$, we have, when $n\to\infty$,
\begin{itemize}
\item[(1)] $a_{X/\F_q,n}\to N_1-1;$
\item[(2)] $b_{X/\F_q,n}\to N_2+N_1^2-2N_1,$
\end{itemize}
where $N_i:=\#X(\F_{q^i})$ denotes the number of $\F_{q^i}$-rational points of $X$.
\end{prop}
\bp
Recall that 
$$
\bc
a_{X/\F_q,n}=\frac{\a_{X/\F_q,n}(n)}{\a_{X/\F_q,n}(0)}-(Q+1);
\\
b_{X/\F_q,n}=(Q-1)\frac{\b_{X/\F_q,n}(0)}{\a_{X/\F_q,n}(0)}+2Q-(Q+1)\frac{\a_{X/\F_q,n}(n)}{\a_{X/\F_q,n}(0)}.
\ec
$$
This means
$$
a_{X/\F_q,n}=\frac{\a_{X/\F_q,n}(n)}{\a_{X/\F_q,n}(0)}-Q-1\to N_1+Q-Q-1=N_1-1.
$$
As for $b_{X/\F_q,n}$, we have 
$$
\ba
b_{X/\F_q,n}&=(Q-1)\frac{\b_{X/\F_q,n}(0)}{\a_{X/\F_q,n}(0)}+2Q-(Q+1)\frac{\a_{X/\F_q,n}(n)}{\a_{X/\F_q,n}(0)}
\\&=Q\frac{\b_{X/\F_q,n}(0)}{\a_{X/\F_q,n}(0)}-\frac{\b_{X/\F_q,n}(0)}{\a_{X/\F_q,n}(0)}+2Q-Q\frac{\a_{X/\F_q,n}(n)}{\a_{X/\F_q,n}(0)}-\frac{\a_{X/\F_q,n}(n)}{\a_{X/\F_q,n}(0)}
\\&\to2(q-1)\whv_1+2q-2N_1
\\&=2(q^2+qa_{X/\F_q}+b_{X/\F_q}+a_{X/\F_q}+1)+2q-2N_1
\\&=N_2+N_1^2-2N_1
\ea
$$
as wanted.
\ep

\begin{thm}[Asymptotic RH for genus 2 curves]
Assume $g=2$, for $n \gg 0$, RH holds for $\zeta_{X/\F_q,n}(s)$.
In particular, if we write $P_{X/\F_q,n}(T)$ as
$$
\ba
\frac{P_{X/\F_q,n}(T)}{\a_{X/\F_q,n}(0)}&=1+a_{X/\F_q,n}T+b_{X/\F_q,n}T^2+a_{X/\F_q,n}QT^3+Q^2T^4
\\&=(1+c_{1,X/\F_q,n}T+QT^2)(1+c_{2,X/\F_q,n}T+QT^2).
\ea
$$
Then we have, for $i\in\{1,2\}$,
$$
\frac{c_{i,X/\F_q,n}}{2\sqrt{Q}}\to \pm\frac{\sqrt{2}}{2},\qquad n\to\infty.
$$
\end{thm}
\bp
According to the estimate of coefficients $a_{X/\F_q,n}$ and $b_{X/\F_q,n}$ of $P_{X/\F_q,n}(T)$, when $n$ goes to infinity, the coefficients will converge to a constant number, namely, $a_{X/\F_q,n}$ goes to $N_1-1$ and $b_{X/\F_q,n}$ goes to $N_2+N_1^2-2N_1$. Moreover, since 
$$
\bc
c_{1,X/\F_q,n}+c_{2,X/\F_q,n}=a_{X/\F_q,n},
\\
c_{1,X/\F_q,n}\times c_{2,X/\F_q,n}=b_{X/\F_q,n}-2Q,
\ec
$$
this means for $i\in\{1,2\}$,
$$
\frac{c_{i,X/\F_q,n}}{2\sqrt{Q}}\to \pm\frac{\sqrt{2}}{2},\qquad n\to\infty.
$$
\ep

\noindent
{\bf Acknowledgement}: The author wishes to thank Professor L. Weng for many helpful discussions and for suggesting the problem studied in this paper. This work was supported by WISE program (MEXT) at Kyushu University.

\vfill
\noindent\text{Zhan SHI}\\
\text{Graduate Program of Mathematics for Innovation}\\
\text{Kyushu University}\\
\text{Fukuoka, Japan}\\
\texttt{shi.zhan.655@s.kyushu-u.ac.jp}
\end{document}